\documentclass{amsart}

\usepackage{amssymb}

\newtheorem{thm}{Theorem}

\newtheorem{lem}[thm]{Lemma}
\newtheorem{cor}[thm]{Corollary}

\newtheorem{prop}[thm]{Proposition}

\newtheorem{conj}[thm]{Conjecture}
   
\theoremstyle{definition}
\newtheorem{defn}[thm]{Definition}

\newtheorem{say}[thm]{}

\newtheorem{ques}[thm]{Question}    

\newtheorem{rem}[thm]{Remark}          

\newtheorem{ack}{Acknowledgments}        
\newtheorem{notation}[thm]{Notation}   
  
\newtheorem{defn-thm}[thm]{Definition--Theorem}  
\newtheorem{defn-lem}[thm]{Definition--Lemma}  

\theoremstyle{remark}

\let \cedilla =\c
\renewcommand{\c}[0]{{\mathbb C}}  

\renewcommand{\o}[0]{{\mathcal O}} 

\renewcommand{\r}[0]{{\mathbb R}} 

\renewcommand{\a}[0]{{\mathbb A}}

\newcommand{\p}[0]{{\mathbb P}}
\newcommand{\f}[0]{{\mathbb F}}
\newcommand{\q}[0]{{\mathbb Q}}
\newcommand{\map}[0]{\dasharrow}
\newcommand{\qtq}[1]{\quad\mbox{#1}\quad}
\newcommand{\spec}[0]{\operatorname{Spec}}

\newcommand{\mult}[0]{\operatorname{mult}}

\newcommand{\supp}[0]{\operatorname{Supp}}

\newcommand{\ex}[0]{\operatorname{Ex}}

\newcommand{\jac}[0]{\operatorname{Jac}}

\newcommand{\simq}[0]{\sim_{\q}}

\def\cprime{$'$}

\begin{document}
\bibliographystyle{amsalpha}

\title{Which powers of holomorphic functions are integrable?}
\author{J\'anos Koll\'ar}

\maketitle
\today

The aim of this lecture is to investigate the following,
rather elementary, problem:

\begin{ques}
Let $f(z_1,\dots,z_n)$ be a holomorphic function on an
open set $U\subset \c^n$. For which $t\in \r$ is
$|f|^t$ locally integrable? 
\end{ques}

The positive values of $t$ pose no problems, for these
$|f|^t$ is even continuous. If $f$ is nowhere zero
on $U$ then again $|f|^t$ is continuous for any $t\in \r$.
Thus the question is only interesting near the
zeros of $f$ and for negative values of $t$.
More generally, if $h$ is an invertible function then
$|f|^t$ locally integrable iff $|fh|^t$ is locally integrable.
Thus the answer to the question depends only on the
hypersurface $(f=0)$ but not on the actual equation.
(A hypersurface $(f=0)$ is not just the
set where $f$ vanishes. One must also remember the  vanishing
multiplicity for
each irreducible component.)

It is traditional to change the question a little 
and work with $s=-t/2$ instead.
Thus we fix a point $p\in U$ and study the values 
$s$ such that $|f|^{-s}$ is $L^2$ in a neighborhood of $p$.
It is not hard to see that there is a
largest value $s_0$  (depending on $f$ and $p$)
such that 
$|f|^{-s}$ is $L^2$ in a neighborhood of $p$
for $s<s_0$ but not $L^2$ 
for $s>s_0$.  Our aim is to study this  ``critical value''
$s_0$. 

\begin{defn} Let $f$ be a holomorphic function in a
neighborhood of a point $p\in \c^n$.
The {\it log canonical threshold} 
or {\it complex singularity exponent} of $f$ at $p$ is
the number $c_p(f)$ such that
\begin{enumerate}
\item[$\bullet$] $|f|^{-s}$ is $L^2$ in a neighborhood of $p$
for $s<c_p(f)$, and
\item[$\bullet$]  $|f|^{-s}$ is not $L^2$ in any neighborhood of $p$
for $s>c_p(f)$.
\end{enumerate}
It is convenient to set $c_p(0)=0$. 
\end{defn}

The name ``log canonical threshold'' comes from algebraic geometry.
I don't know who studied these numbers first. The concept is probably too
natural to have a well defined inventor. 
It appears in the works of Schwartz, H\"ormander,
{\L}ojasiewicz and Gel{\cprime}fand as
 the ``division problem for distributions''; 
see \cite{sch, hor, loj, gel}.
  The general question is considered
by Atiyah \cite{ati} and Bernstein \cite{ber}.
The connections with singularity theory were explored
by the Arnol{\cprime}d school  and summarized in \cite{agv}. 
See \cite{k-pairs} for another 
survey and for further connections.

Algebraic geometers became very much interested in 
log canonical thresholds when Shokurov \cite{shok88}
 discovered that some subtle
properties of log canonical thresholds, especially his conjecture
 (\ref{acc.conj}), are
connected with the general MMP   (=minimal model program).
These connections were systematized and further developed  
in \cite[Secs.17--18]{f&a}. In this general framework,
considering only smooth complex spaces is not natural.
In fact, the inductive theorems require the consideration
of cases when $f$ is holomorphic on a singular
complex space. 
At the end, the singular versions of the 
ACC (=ascending chain condition) Conjecture
(\ref{acc.conj}) and the Accumulation Conjecture (\ref{accum.conj})
emerged as the main open problems.

A novel approach to log canonical thresholds 
on manifolds was proposed by
de~Fernex and Musta{\cedilla{t}}{\v{a}}
 in \cite{df-m}.  
 They rely on 
non-standard methods (ultraproducts etc.) and the formula for
log canonical thresholds using arc-spaces  \cite{mus-jet}.
The end result is the proof of the smooth version of the
Accumulation Conjecture 
for decreasing limits.

The aim of this lecture is three-fold.
First, I give an elementary introduction to
log canonical thresholds.
The second part is a presentation of
the  proof in \cite{df-m}  using ``traditional'' methods and the
original definition of log canonical thresholds relying on
discrepancies of divisors as in (\ref{compute.lcth}.4).
At its heart, however, the proof in (\ref{limit.series.thm})
  is the same as in \cite{df-m}.
Third, I show how to use the existence of minimal models
\cite{bchm} to establish a  part of the ACC conjecture.
This in turn is enough to complete the proof of the smooth version of
the Accumulation Conjecture. 

There are three, quite distinct, approaches to 
log canonical thresholds.
\begin{enumerate}
\item[$\bullet$] Study the relationship of $c_p(f)$
and the singularity $p\in (f=0)$. 
\item[$\bullet$]  Study the function $f\mapsto c_0(f)$ 
on the space of all holomorphic functions.
\item[$\bullet$]  Study the set of all possible values
$c_0(f)\in \r$.
\end{enumerate}

The log canonical threshold is related to  other invariants
of singularities in many ways; see \cite{k-pairs} for a survey.
However, I will not say anything about these here, mainly because
in higher dimensions  these connections have not yet proved useful in the
study of the other two problems.

These notes start at an elementary level. I tried hard to avoid
the algebraic methods and terminology. However, starting with Section 6,
I switch to the language of divisors since it is  better suited
to handle the general singular case.

\section{Main conjectures}

It was Shokurov in \cite{shok88} who first proposed to look at all possible
values of log canonical thresholds in a fixed dimension and suggested that
these sets, though rather complicated, have 
remarkable properties. The original questions were extended and
further developed in \cite[Sec.18]{f&a}.

\begin{defn}\label{c0.defn} Let ${\mathcal {HT}}_n$ be the set of 
log canonical thresholds of all possible $n$-variable holomorphic functions.
That is, 
$$
{\mathcal{HT}}_n:=\Bigl\{ c_0(f): f\in \o_{0,\c^n}\Bigr\}.
$$
The notation suggests that we are talking about
{\it h}ypersurface {\it t}hresholds.
As we see in Section \ref{taylor.sec}, we get the same set if instead we let
$f$ run through all polynomials or all formal power series
over any algebraically closed field of characteristic 0.
As far as I know, the answer could be the same
if we look at polynomials over any field (e.g.\ $\q$ or even $\f_p$).

The sets  ${\mathcal {HT}}_n$ are different from the sets
${\mathcal T}_n$ used in \cite{k-pairs} and in \cite{df-m}
(which are also different from each other).

The paper \cite{df-m} considers log canonical thresholds when a
single holomorphic function $f$ is replaced by
$\max\bigl\{|f_1|,\dots, |f_r|\bigr\}$
where the $f_i$ are holomorphic. It is easy to rework the
results of this note in their more general setting.

The lectures \cite{k-pairs} consider log canonical thresholds
for functions on singular complex spaces.
The present methods apply to that case if there is an
a priori bound on the appearing singularities. 
For instance, the proofs work if we assume that $X$ has 
only hypersurface singularities.
Using the standard covering and partial resolutions tricks
(for instance, as in \cite[Sec.5]{lsgt}),
this is a manageable  limitation in dimensions $\leq 3$.
Unfortunately, this is a rather unnatural restricion
in connection with the higher dimensional MMP.

Note that  $|z|^{-s}$ is $L^2$ iff $s<1$. From this we conclude
that, for a  1-variable holomorphic function $f(z)$,
$$
c_p\bigl(f(z)\bigr)=\frac1{\mult_pf}.
$$
In particular, 
$$
{\mathcal{HT}}_1=\bigl\{1, \tfrac12, \tfrac13, \dots, 0\bigr\}.
\eqno{(\ref{c0.defn}.1)}
$$
The 2-variable case is already quite subtle, but
 we see in (\ref{newton.est.thm}.5) that
$$
{\mathcal{HT}}_2=\Bigl\{\frac{c_1+c_2}{c_1c_2+a_1c_2+a_2c_1}:
 a_i+c_i\geq \max\{2,a_{3-i}\}\Bigr\}\cup\{0\}.
\eqno{(\ref{c0.defn}.2)}
$$
\end{defn}

Although the sets ${\mathcal{HT}}_n$ are not known for $n\geq 3$,
and a complete listing as in (\ref{c0.defn}.2)
may not even be interesting, they are conjectured to 
have remarkable properties.
 The following are the basic results and questions:

\begin{prop}\label{rtl.01.prop}
 All log canonical thresholds are rational and lie between $0$ and $1$.
That is, 
$$
{\mathcal{HT}}_n\subset \q\cap [0,1].
$$
\end{prop}

This is proved in (\ref{simple.est.lem}) and
(\ref{compute.lcth}).

The key question in this area is the following:

\begin{conj}[ACC conjecture, smooth version]\label{acc.conj} 
For any $n$ there is no infinite increasing subsequence
in ${\mathcal{HT}}_n$.
\end{conj}

Note by contrast, that by (\ref{dsum.prop}.1), any 
rational number between 0 and 1
is the log canonical threshold of some function for some $n$.

There are many decreasing sequences of log canonical thresholds, and
the following  conjecture 
\cite[8.21.2]{k-pairs} describes their limit points:

\begin{conj}[Accumulation conjecture, smooth version] \label{accum.conj}
The set of accumulation points of ${\mathcal{HT}}_n$
is ${\mathcal{HT}}_{n-1}\setminus\{1\}$.
\end{conj}

It is easy to see  (\ref{dsum.prop}.2) that the 
accumulation points of ${\mathcal{HT}}_n$
contain ${\mathcal{HT}}_{n-1}\setminus\{1\}$.
The main result of this note is to prove
that Conjecture \ref{accum.conj} almost holds:

\begin{thm} \label{main.accu.thm}
The set of accumulation points of ${\mathcal{HT}}_n$
is either  ${\mathcal{HT}}_{n-1}\setminus\{1\}$ or ${\mathcal{HT}}_{n-1}$.
\end{thm}

By (\ref{rtl.01.prop}), $1\in \r$ can not be a limit of
a decreasing sequence of log canonical thresholds.
As a special case of the ACC conjecture,
$1$ can not be a limit of
an increasing sequence of log canonical thresholds either.
Equivalently, for a fixed dimension, no log canonical threshold lies in
an interval  $(1-\epsilon_n,1)$ for some $\epsilon_n>0$. I call this 
 special case the Gap conjecture \cite[8.16]{k-pairs}.

\begin{conj}[Gap conjecture, smooth version]  \label{gap.conj}
For every $n$ there is an  $\epsilon_n>0$ such that
for any $f(z_1,\dots,z_n)$ that is holomorphic on
the closed unit ball $B$,
$$
\int_B \frac1{(f\bar f)^{1-\epsilon_n}}\ dV<\infty 
\quad\Rightarrow\quad 
\int_B \frac1{(f\bar f)^{1-\epsilon}}\ dV<\infty 
\quad \forall\  0<\epsilon<\epsilon_n.
$$
\end{conj}

Various forms of the Gap conjecture 
are important in the construction of Einstein metrics
as in \cite{bgk}. Ultimately, a gap conjecture type result lies behind
the stabilization theorems in \cite{k-em1, k-em2}.

As noted in \cite{df-m}, the Gap conjecture and
(\ref{main.accu.thm}) imply the ACC conjecture.
However, to obtain the ACC conjecture in a fixed dimension,
one needs the Gap conjecture in all dimensions.

There is even a conjecture about the precise value of the optimal
$\epsilon_n$.

Consider the sequence defined recursively by 
$c_{k+1}=c_1\cdots c_k+1$
starting with $c_1=2$.(It is called
Euclid's or Sylvester's sequence,
see \cite[Sec.4.3]{gkp} or \cite[A000058]{sequences}.)
It starts as
$$
2,3,7,43,1807, 3263443, 10650056950807,...
$$
It is easy to see  that 
$$
\sum_{i=1}^n\frac1{c_i}=1-\frac1{c_{n+1}-1}= 1-\frac{1}{c_1\cdots c_n}.
$$
In particular, by (\ref{dsum.prop}.1),
$$
c_0\Bigl(z_1^{c_1}+\cdots+z_n^{c_n}\Bigr)=1-\frac{1}{c_{n+1}-1}.
$$
It is conjectured that this is the worst example, that is,
the optimal value for $\epsilon_n$ in (\ref{gap.conj}) is
$$
\epsilon_n=\frac{1}{c_{n+1}-1}.
\eqno{(\ref{gap.conj}.1)}
$$

\begin{say}[Known special cases] 

As we noted, 
${\mathcal{HT}}_1$ and
${\mathcal{HT}}_2$ are known. From these one can read off
all the above conjectures for $n\leq 2$.
In particular, we get that
$\epsilon_1=\frac{1}{2},\ \epsilon_2=\frac{1}{6}$.
 The value
$\epsilon_3=\frac{1}{42}$ is computed in \cite[5.5.7]{lsgt},
 essentially through a classification of the possible normal forms
of singularities with log canonical threshold near 1.

The set ${\mathcal{HT}}_3$ is still not known,
but \cite{kuwata} determined
${\mathcal{HT}}_3\cap[\frac56,1]$
and \cite{prok} computed all accumulation points of
${\mathcal{HT}}_3$ that lie in $[\frac12,1]$.

The ACC conjecture for ${\mathcal{HT}}_3$ was proved by
\cite{alexeev} and the Accumulation conjecture
by \cite{mck-pr}. Both of these papers deal with the general
singular case and rely heavily on the MMP in dimension 3.
The relevant parts of the MMP are now known in all dimensions \cite{bchm}.
A missing ingredient in higher dimensions is the
Alexeev-Borisov-Borisov conjecture \cite{al94}. 
Even stating it would lead us quite far.
The toric cases are treated in \cite{bo-bo}.

 \cite{sound} proved that, with $\epsilon_n$ as in (\ref{gap.conj}.1),  
$c_0(z_1^{a_1}+\cdots+z_n^{a_n})$ can not lie in
$(1-\epsilon_n, 1)$ for any
$a_1,\dots,a_n$. That is, if
$$
\frac1{a_1}+\cdots+\frac1{a_n}<1 \qtq{then}
 \frac1{a_1}+\cdots+\frac1{a_n}\leq 1-\epsilon_n.
$$
\end{say}

\section{Computing and estimating $c_0(f)$}

In this section we discuss how to determine or bound the
log canonical threshold.  The basic result (\ref{compute.lcth}),
first observed by Atiyah, gives a formula for
$c_0(f)$ in terms of an embedded resolution of the
hypersurface $(f=0)$.  It is not easy to construct
embedded resolutions, but even simple-minded partial
resolutions frequently 
give good upper bounds for $c_0(f)$. Estimates using the
Newton polygon are especially easy to obtain and to use.
It is much harder to get good lower bounds. 

\begin{lem}\label{simple.est.lem} If $f(p)\neq 0$ then
$c_p(f)=+\infty$. If $f(p)= 0$ then $0\leq c_p(f)\leq 1$. 
\end{lem}

Proof. The first claim is clear. Thus assume that  $f(p)= 0$.
As we noted in (\ref{c0.defn}), for a  1-variable holomorphic function $f(z)$
we have 
$c_p\bigl(f(z)\bigr)=\frac1{\mult_pf}$.

In the several variable case, pick a smooth point $p'$ near $p$ on the
hypersurface $(f=0)$. We can choose local coordinates near $p'$
such that $f=(\mbox{unit})z_1^m$ for some $m$. By Fubini this
shows that  $c_p(f)\leq c_{p'}(f)=1/m$. \qed
\medskip

As we see in (\ref{compute.lcth}.5) and (\ref{ot-inv}.1),
in several variables one can only get 
inequalities relating $c_0(f)$ and the
multiplicity:
$$
\frac1{\mult_pf}\leq c_p\bigl(f(z_1,\dots,z_n)\bigr)\leq \frac{n}{\mult_pf}.
$$

\begin{say}[Computing the log canonical threshold]\cite{ati}
\label{compute.lcth}
Set $\omega=dz_1\wedge\cdots\wedge dz_n$. 
Then $|f|^{-s}$ is locally $L^2$ iff, on any compact
$K\subset U$, the integral
$$
\int_K (f\bar f)^{-s} \ \omega\wedge \bar{\omega} \qtq{is finite.}
\eqno{(\ref{compute.lcth}.1)}
$$
(We can ignore the power of $\sqrt{-1}$ that makes this integral real.)
Let $\pi:X\to U$ be a proper bimeromorphic morphism.
We can  rewrite the above integral as
$$
\int_K (f\bar f)^{-s} \ \omega\wedge \bar{\omega}=
\int_{\pi^{-1}(K)} \bigl((f\circ \pi)\overline{(f\circ \pi)}\bigr)^{-s}
\  \pi^*\omega\wedge 
\pi^*\bar{\omega}.
\eqno{(\ref{compute.lcth}.2)}
$$
The aim  now is to choose $\pi$ such that
the local structure of $f\circ \pi$ and of $\pi^*\omega$ becomes simple.
The best one can do is to take an embedded  resolution of singularities
for $(f=0)$. This is a
proper bimeromorphic morphism $\pi:X\to U$
such that $X$ is a smooth complex manifold and
the zero set of $f\circ \pi$  plus the exceptional set of $\pi$
is a normal crossing divisor.
That is,
at any point $q\in X$ we can choose local coordinates
$x_1,\dots, x_n$ such that
$$
 f\circ \pi=(\mbox{invertible})\prod_i x_i^{a(i,q)}
\qtq{and}
\pi^*\omega=(\mbox{invertible})\prod_i x_i^{e(i,q)}\cdot
dx_1\wedge\cdots\wedge dx_n,
$$
where $a(i,q)=\mult_{(x_i=0)}(f\circ \pi)$ and
$e(i,q)=\mult_{(x_i=0)}\jac \pi$.
Here $\jac$ denotes the complex Jacobian
$$
\jac\pi=\det\Bigl( \frac{\partial z_i}{\partial x_j}\Bigr).
$$
Thus the integral (\ref{compute.lcth}.2) is finite near $q\in X$ iff
$$
\int\cdots\int \prod_i (x_i\bar{x_i})^{e(i,q)-s\cdot a(i,q)}\ dV \ =
 \pm \prod_i\ \int(x_i\bar{x_i})^{e(i,q)-s\cdot a(i,q)}\ dx_i\wedge d\bar x_i
\eqno{(\ref{compute.lcth}.3)}
$$
is finite. 
This holds iff $e(i,q)-s\cdot a(i,q)>-1$ for every $i$, that is,
when $s< (e(i,q)+1)/a(i,q)$. 
This gives the formula for the  log canonical threshold:
$$
c_p(f)=\min\left\{ \frac{1+\mult_{E}\jac \pi}{\mult_{E} (f\circ \pi)}:
  \ \mbox{ for those $E$ such that }
\ p\in \pi(E)\right\}.
\eqno{(\ref{compute.lcth}.4)}
$$
In principle we can take the minimum over all divisors
$E\subset X$ such that $p\in \pi(E)$. 
However, only the exceptional divisors of $\pi$ and the
(birational transforms of) irreducible components of $(f=0)$
are interesting. For all other $E$, 
$\mult_{E}\jac \pi=\mult_{E} (f\circ \pi)=0$
and their contribution to (\ref{compute.lcth}.4)
is $1/0=+\infty$.

It is customary to view $\mult_E(f\circ \pi)$ as a valuation on functions
on $\c^n$ and drop $\pi$ from the notation. Thus we write
$$
\mult_E f\qtq{instead of} \mult_E(f\circ \pi).
$$

First of all,  the formula (\ref{compute.lcth}.4) shows that
the log canonical threshold of $f$ is always a rational number,
completing the proof of (\ref{rtl.01.prop}). 

Second, it gives us ways to compute or at least estimate
$c_0(f)$. 

In  many cases it is not hard to guess which exceptional divisor
computes the log canonical threshold 
(that is, achieves equality in (\ref{compute.lcth}.4)), and to write down a  
bimeromorphic morphism $\pi:X\to U$
where this divisor appears. 
This way we can get upper bounds for $c_0(f)$. 
Note that we do not need  to arrange
that $X$ be smooth or that $\pi$ be proper. 
Any bimeromorphic morphism $\pi_1:X_1\to U$ can be 
completed to a proper bimeromorphic morphism $\pi_2:X_2\to U$ and
then, by resolution of singularities,  
to a proper bimeromorphic morphism $\pi_3:X_3\to U$
which is an embedded resolution as required for (\ref{compute.lcth}.4).

For instance, let $\pi:B_0\c^n\to \c^n$ be the blow up of the origin
with exceptional divisor $E\cong \p^{n-1}$.
Then $\mult_{E}\jac \pi=n-1$ and $\mult_{E} (f\circ \pi)=\mult_0f$.
This gives the simple estimate
$$
c_0\bigl(f(z_1,\dots,z_n)\bigr)\leq \frac{n}{\mult_0f}.
\eqno{(\ref{compute.lcth}.5)}
$$
\end{say}

\begin{say}[Formal power series] \label{formal.defn}
The formula (\ref{compute.lcth}.4) makes it possible
to define the log canonical threshold for a formal
power series $f\in k[[z_1,\dots,z_n]]$ over any
field $k$ of characteristic 0. Indeed, resolution of
singularities is known for complete local rings  \cite{tem}
and then (\ref{compute.lcth}.4) makes sense. 
It is easy to see that the resulting $c_0(f)$
is independent of the resolution.

An alternative definition of $c_0(f)$ using arc spaces
is given in \cite{mus-jet}.
\end{say}

An especially convenient estimate is obtained
 using the  Newton polygon.

\begin{defn}[Newton polygon]\label{newton.polyg} 
Let  $F=\sum a_{I}{\mathbf x}^I$ be a polynomial or power series 
  in $n$-variables.
The  {\it Newton polygon}
 of $f$ (in the chosen coordinates $x_1,\dots,x_n$)
is obtained as follows.

In $\r^n$  we mark  the point $I=(i_1,\dots,i_n)$ with a big dot if
$a_{I}\neq 0$.
Any other monomial $x^{I'}$ with $I'\geq I$ coordinatewise will not be of
``lowest order'' in any sense, so we also mark these.
(In the  figure below these markings are invisible.)

The  {\it Newton polygon}
is the boundary of the convex hull of the resulting infinite set of
marked  points.
$$
\begin{picture}(100,80)(0,20)
\put(0,0){\line(1,0){80}}
\put(0,0){\line(0,1){80}}

\put(0,0){\circle*{2}}
\put(0,10){\circle*{2}}
\put(0,20){\circle*{2}}
\put(0,30){\circle*{2}}
\put(0,40){\circle*{2}}
\put(0,50){\circle*{2}}
\put(0,60){\circle*{2}}
\put(0,70){\circle*{2}}

\put(10,0){\circle*{2}}
\put(20,0){\circle*{2}}
\put(30,0){\circle*{2}}
\put(40,0){\circle*{2}}
\put(50,0){\circle*{2}}
\put(60,0){\circle*{2}}
\put(70,0){\circle*{2}}
\put(80,0){\circle*{2}}

\put(0,70){\circle*{4}}
\put(20,30){\circle*{4}}
\put(40,10){\circle*{4}}
\put(50,30){\circle*{4}}
\put(60,0){\circle*{4}}

\put(0,70){\line(1,-2){20}}
\put(20,30){\line(1,-1){20}}
\put(40,10){\line(2,-1){20}}
\end{picture}
\begin{array}{c}
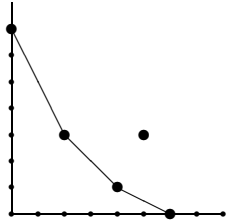
\mbox{The Newton polygon of}\\
 y^7+y^3x^2+y^3x^5+yx^4+x^6.\\
{\ }\\
{\  }\\
{\ }
\end{array}
$$
\end{defn}

\begin{say}[Estimating $c_0(f)$ using ${\rm Newton}(f)$]

Let $\sum_i a_ix_i=d$ be the equation of a face
of ${\rm Newton}(f)$. We can assume that the $a_i$ are relatively prime
positive integers.

Then $a_1,\dots,a_n$ can be the first column of an $n\times n$
invertible integral matrix $M=(a_{ij})$. 
Consider the map
$$
\pi: \c_{x_1}\times (\c^*)^{n-1}_{x_2,\dots,x_n} \to \ \c^n_{z_1,\dots,z_n}
\qtq{given by} z_i=\textstyle{\prod}_j x_j^{a_{ij}}.
$$
The inverse of $M$ defines the inverse of $\pi$
on the open subset $(\c^*)^{n}\subset \c^n_{z_1,\dots,z_n}$.
We concentrate on the exceptional divisor $E:=(x_1=0)$. 

Note that $\pi^*dz_i=
\bigl(\prod_j x_j^{a_{ij}}\bigr)\bigl(\sum_ja_{ij}\frac{dx_j}{x_j}\bigr)$,
hence
$$
\pi^*\bigl(dz_1\wedge\cdots\wedge dz_n\bigr)=
\bigl(\textstyle{\prod}_{ij}x_j^{a_{ij}}\bigr)\cdot \det M
\cdot (x_1\cdots x_n)^{-1}\cdot (dx_1\wedge\cdots\wedge dx_n\bigr).
$$
Thus the Jacobian of $\pi$ vanishes along
$E$ with multiplicity  $-1+\sum_ia_{i1}=-1+\sum_i a_i$.

If $\prod_i z_i^{b_i}$ is any monomial occurring in $f$
then $\sum_i a_ib_i\geq d$
since $(b_1,\dots,b_n)$ lies above $\sum_i a_ix_i=d$. On the other hand
$$
\bigl(\textstyle{\prod}_i z_i^{b_i}\bigr)\circ \pi=
\textstyle{\prod}_j x_j^{A_j}\qtq{where}
A_j=\sum_i b_ia_{ij}
$$
and so it vanishes along $E$ with 
multiplicity  $A_1=\sum_i b_ia_{i1}=\sum_i b_ia_{i}\geq d$.
Thus we conclude that 
$$
c_0(f)\leq \frac{\sum_i a_i}{d}.
$$

\end{say}

This inequality is equivalent to the first part of
the next theorem. For the proof of the second part  see \cite{kous}
and for the third 
 \cite{varch} or \cite[6.40]{ksc}.

\begin{thm} \label{newton.est.thm}
Let $f(z_1,\dots,z_n)$ be a holomorphic function near the origin $0\in \c^n$.
Let ${\rm Newton}(f)$ be the Newton polygon of $f$.
\begin{enumerate}
\item The vector  $(1/c_0(f), \cdots, 1/c_0(f))$ is on or above 
${\rm Newton}(f)$.
\item Fix ${\rm Newton}(f)$ and assume that the coefficients of 
the monomials in $f$ are general. Then 
 $(1/c_0(f), \cdots, 1/c_0(f))$ is on 
${\rm Newton}(f)$.
\item  If $n=2$ then one can choose  local analytic
coordinates $(x,y)$
such that $(1/c_0(f),  1/c_0(f))$ is on ${\rm Newton}(f)$.
\end{enumerate}
\end{thm}

This gives an easy way to construct many different
log canonical thresholds. Take $m\leq n$ linearly
 independent nonnegative vectors
${\mathbf a}_i=(a_{ij})$ such that their convex hull 
contains a vector of the form
$(1/c, \cdots, 1/c)$. Then, for general $b_i\in \c$,
$$
c_0\Bigl(\textstyle{\sum}_i \ b_i \textstyle{\prod}_j x_j^{a_{ij}}\Bigr)=c.
\eqno{(\ref{newton.est.thm}.4)}
$$
After a suitable change of coordinates we can even
assume that all $b_i=1$.
These are the log canonical thresholds that can be computed
as in (\ref{compute.lcth}) using a resolution $X\to \c^n$ which is
{\it toric}, that is, equivariant with respect to the standard
$(\c^*)^n$-action on $\c^n$. 
The values produced by (\ref{newton.est.thm}.4) give a large subset
of ${\mathcal{HT}}_n$. It is possible that in fact these values give all of 
${\mathcal{HT}}_n$, but 
I know of no reasons why this should be true.
Note, however, that (\ref{newton.est.thm}.3) definitely fails
already for $n=3$; see \cite[6.45]{ksc} for an example.

Computing the case when $(1/c_0(f),  1/c_0(f))$ is on 
the edge of the Newton polygon between the points
$(a_1, a_2+c_2)$ and $(a_1+c_1, a_2)$ shows that
any 2-variable log canonical threshold can be written as
$$
c_0(f)=
\frac{c_1+c_2}{c_1c_2+a_1c_2+a_2c_1}
\eqno{(\ref{newton.est.thm}.5)}
$$
where $a_1+c_1\geq \max\{2,a_2\}$  and 
$a_2+c_2\geq \max\{2,a_1\}$, or it is $0$ or $1$.

\section{Basic properties}\label{basic.props.sec}

In this section we collect the known important
properties of log canonical thresholds.
I state everything for formal power series, as needed
for our proofs. See (\ref{formal.say}) for comments on the proofs
in this setting.

The next result is proved for holomorphic functions with
isolated critical points in \cite[II.13.3.5]{agv}.
The proof given in \cite[8.21]{k-pairs} works in general.

\begin{prop} \label{dsum.prop}
Let $f({\mathbf x})$ and
$g({\mathbf y})$ be power series in disjoint sets of variables.
Then
$$
c_0\bigl(f({\mathbf x})\oplus g({\mathbf y})\bigr)=\min\bigl\{ 1, 
c_0(f)+c_0(g)\bigr\},
$$
where $\oplus$ denotes the sum in disjoint sets of variables.
\end{prop}

As a  corollary we obtain that
$$
c_0\bigl(z_1^{a_1}+\cdots+z_n^{a_n}\bigr)=
\min\Bigl\{ 1, \frac1{a_1}+\cdots +\frac1{a_n}\Bigr\},
\eqno{(\ref{dsum.prop}.1)}
$$
and if $c_0\bigl(f({\mathbf x})\bigr)<1$ and $m\gg 1$ then
$$
c_0\bigl(f({\mathbf x})+y^m\bigr)=c_0\bigl(f({\mathbf x})\bigr)+\tfrac1{m}.
\eqno{(\ref{dsum.prop}.2)}
$$

A simple but important property is that the
 log canonical threshold  gives a metric on the space of
power series that vanish at the origin, 
\cite[Thm.2.9]{dem-ko} or \cite[8.19]{k-pairs}.

\begin{thm}\label{dk.thm}
 Let $f({\mathbf x}), g({\mathbf x})\in k[[x_1,\dots,x_n]]$
 be power series.
Then
$$
c_0(f+ g)\leq  c_0(f)+c_0(g).
$$
\end{thm}

Applying this to the Taylor polynomials $t_m(f)$ of $f$ and
to $f-t_m(f)$, and using
(\ref{compute.lcth}.5),  we get the following  uniform approximation result.

\begin{cor} \label{taylor.app.cor}
 Let $f({\mathbf x})\in k[[x_1,\dots,x_n]]$
 be a power series and $t_m(f)$ its degree $m$ Taylor polynomial. Then
$$
\Bigl|c_0\bigl(f\bigr)-c_0\bigl(t_m(f)\bigr)\Bigr|\leq \frac{n}{m+1}. \qed
$$
\end{cor}

\begin{notation} In order to indicate the change form the complex
to the algebraic case, I replace $\c^n$ with the affine $n$-space
$\a^n$ defined over some field $k$. Its completion at the origin
is denoted by $\hat\a^n$. It is also
$\spec_k k[[x_1,\dots,x_n]]$.
\end{notation}

The next result is known as the Ohsawa-Takegoshi
extension theorem in complex analysis \cite{oh-ta} and
as the (weak version of) inversion of adjunction
in algebraic geometry
\cite[Sec.17]{f&a}.

\begin{thm} \label{ot-inv} Let $f({\mathbf x})\in k[[x_1,\dots,x_n]]$
 be a power series and $L\subset \hat{\a}^n$ a smooth
submanifold.
Then
$$
c_0\bigl(f|_L)\leq c_0(f).
$$
\end{thm}

For instance, if $L\subset \c^n$ is a general line through the
origin then $\mult_0\bigl(f|_L\bigr)=\mult_0f$
and so we conclude that
$$
c_0(f)\geq \frac1{\mult_0f}.
\eqno{(\ref{ot-inv}.1)}
$$

\begin{say}[Thresholds in families]\label{thr.family.say}
Let $f_x:=\sum_I a_I(x){\mathbf z}^I$ be polynomials in 
${\mathbf z}$ whose coefficients $a_I(x)$
are rational functions on an algebraic variety $X$.
What can we say about the log canonical thresholds $c_0(f_x)$
as a function of $x\in X$?

Pick a generic point $x_g\in X$ and 
take a resolution $\pi_{x_g}:Y_{x_g}\to \a^n$.
Since $\pi_{x_g},Y_{x_g}$ and the exceptional divisors $E^i_{x_g}$
 are defined over the generic point of $X$,
 there is a 
Zariski open subset  $X^0\subset X$
such that, by specialization,  for every $x\in X^0$
we obtain $\pi_{x}:Y_{x}\to \a^n$ with  exceptional divisors $E^i_{x}$.
Moreover, we may assume that, 
for every $i$, the multiplicities of
$\jac \pi_x$ and of $f_x\circ \pi_x$ along $E^i_x$ do not depend on $x\in X^0$.
In particular, $c_0(f_x)$ is also independent of
$x\in X^0$.

Repeating the argument with $X$ replaced by $X\setminus X^0$,
we conclude that $c_0(f_x)$ is a constructible function of
$x\in X$. That is, its level sets are finite unions of
locally closed subvarities.

It is also easy to see that $c_0(f_x)$ is a 
lower semi continuous function of $x\in X$, cf.\ \cite{varch}.
\end{say}

In the complex analytic case, a more precise version
of lower semi continuity
 is proved in \cite[0.2]{dem-ko}:

\begin{thm} Assume that $f_t({\mathbf z})$ converges uniformly to 
$F({\mathbf z})$
in a compact neighborhood $B$ of $0$. Fix $s<c_0(F)$. Then
$$
\frac1{|f_t|^s}\qtq{converges to} \frac1{|F|^s}\qtq{in $L^2(B)$.}
$$
\end{thm}

\begin{say}[Comments on the formal power series case]
\label{formal.say}

All these results were proved in the algebraic and analytic settings.
The methods are either analytic, or, as the proofs of 
inversion of adjunction in \cite[Sec.17]{f&a}
and \cite{kawakita},  rely ultimately on a relative version of
the Kodaira vanishing theorem. This vanishing is known for
birational maps between varieties and for bimeromorphic maps
between complex spaces. Unfortunately, we would need it
in case the base is  a formal power series ring. 

While the result is no doubt true in this case, the usual proofs
of the Kodaira-type vanishing theorems rely on some
topological/analytic arguments. Thus, a genuinely new proof
may be needed. 

Here I go around this difficulty by a  reduction to the
algebraic case, see
Section \ref{techn.sec}.
This, however, should be viewed as  but a
temporary patch. It is  high time to work out
the whole MMP over an arbitrary base scheme, especially
over complete local rings.

Note that the formal versions of 
(\ref{ot-inv}) and (\ref{dk.thm})
both follow from the algebraic case once we know that
for any formal power series $f$, the log canonical thresholds of its Taylor
approximations converge to $c_0(f)$. That is, if
$$
\lim_{m\to\infty} c_0\bigl(t_m(f)\bigr)=c_0(f)
\eqno{(\ref{formal.say}.1)}
$$
The argument in (\ref{upper.anytail.say}) easily yields the inequality
$$
\limsup_{m\to\infty}\ c_0\bigl(t_m(f)\bigr)\leq c_0(f),
$$
but the other direction relies on inversion of adjunction
(in larger dimensions), creating a vicious circle.

The first complete proof of
(\ref{formal.say}.1)  is  in \cite[2.5]{df-m}
using arc-space techniques. 
\end{say} 

\begin{say}[Proof of (\ref{dsum.prop}) $\wedge$
(\ref{ot-inv}) $\Rightarrow$ (\ref{dk.thm})]\label{dk.proof}

Create disjoint sets of variables for $f({\mathbf x})$ and $g({\mathbf y})$.
Then, by (\ref{dsum.prop}), 
$$
c_0\bigl(f\oplus g\bigr)\leq  c_0(f)+c_0(g).
$$
Note that $f({\mathbf x})+g({\mathbf x})$ is naturally isomorphic to
$f({\mathbf x})\oplus g({\mathbf y})$
restricted to the diagonal
$L:=(x_1-y_1=\cdots=x_n-y_n=0)$.
Thus, by (\ref{ot-inv}),
$$
c_0(f+ g)\leq c_0\bigl(f\oplus g\bigr)\leq  c_0(f)+c_0(g).\qed
$$
\end{say}

\section{Generic limits of power series}

\begin{say} Consider a sequence of holomorphic functions
$f_i$ defined in a neighborhood of  $0\in \c^n$.
Assume that the sequence of log canonical thresholds converges to 
a limit  $c:=\lim_i c_0(f_i)$. Can we write down a 
holomorphic function $f$ such that $c_0(f)=c$
and, in some sense, $f$ is the limit of the functions
$f_i$?

At first sight the answer is no. Even in some very simple cases
when the $f_i$ do converge to a limit, the log canonical threshold
 usually jumps.
For instance, take $f_i(z)=z^2+\frac1{i}z$ and $f(z)=z^2$.
Then $f_i\to f$ uniformly on any compact set, yet
$c_0(f_i)=1$ and $c_0(f)=\frac12$. 

We get a different insight from the log canonical threshold formula using
exceptional divisors (\ref{compute.lcth}.4). Let $\pi:X\to \c^n$ be a 
bimeromorphic map and $E\subset X$  a divisor such that
$\pi(E)=0$. Choose local coordinates 
$x_1,\dots, x_n$ at a general point of $E$
such that $E=(x_1=0)$. If $\pi_i$ are the coordinate functions
of $\pi$ then,
expanding $f\circ \pi$  by powers of ${\mathbf x}$, write
$$
f\bigl(\pi_1({\mathbf x}),\dots, \pi_n({\mathbf x})\bigr)=
\sum_{I} P_I({\mathbf a}, {\mathbf b}){\mathbf x}^I,
$$
where the $P_I$ are polynomials, the ${\mathbf a}$ are
the coefficients of $f$ and the ${\mathbf b}$
are the coefficients of $\pi$.

Note that   $f\circ \pi$ vanishes along
$E$ with  multiplicity $m$ iff
 $x_1^m$ divides $f\circ \pi$.
Equivalently, when $P_I({\mathbf a}, {\mathbf b})=0$
whenever the first coordinate of $I=(i_1,\dots,i_n)$
is less than $m$.

This suggests that we should focus
on  the polynomial relations
between the coefficients $a_{J}$.
This is a key idea that 
de~Fernex and Musta{\cedilla{t}}{\v{a}}  
use to study limits of log canonical thresholds.

An interesting feature of the proof is that
even if we start with a sequence of functions $f_i$
that are  holomorphic on a fixed open set $U$, their
limit is only a formal power series $F(z_1,\dots,z_n)$.
Furthermore, the construction naturally
yields a power series $F$ whose coefficients are not in $\c$
but in  an algebraically closed field $K$ of
countably infinite transcendence degree over $\c$.

Any such field $K$ is isomorphic to $\c$, so at the end
we can replace $F$ with a 
formal power series $F^*(z_1,\dots,z_n)\in \c[[z_1,\dots,z_n]]$
and, using (\ref{anytail.prop}), even  with a polynomial
$P(z_1,\dots,z_n)\in \c[z_1,\dots,z_n]$, but
these steps are rather artificial from the point of view of the proof. 
It is more natural to work 
with formal power series
over an arbitrary field $k$.
\end{say}

\begin{say}[Generic power series] 
 Let $k$ be a field and $k[[ x_1,\dots,x_n]]$
the ring of power series with coefficients in $k$. 
We can view $k[[ x_1,\dots,x_n]]$ as an infinite dimensional
affine space $\a^{\infty}$ over $k$. 
Thus if $f_i(x_1,\dots,x_n)\in k[[ x_1,\dots,x_n]]$
are power series, then we get points $[f_i]\in \a^{\infty}_k$.
Assume now that there is a power series  $F\in K[[ x_1,\dots,x_n]]$
 over a possibly larger field $K$ such that $[F]$ is a 
``generic point'' 
of the ``Zariski closure'' 
$Z\subset\a^{\infty}$ of
the  set  $\bigl\{[f_i]: i\in I\bigr\}$.
The first main result of \cite{df-m} says, roughly, that
$$
c_0(F)=\lim_{j\to\infty}
 c_0(f_{i_j})\qtq{for some subsequence $i_1<i_2<\cdots$.}
$$
One needs to be rather careful with ``Zariski closure'' and
``generic point'' in an infinite dimensional setting.

The non-standard method in \cite{df-m} is used to get
a correct ``generic point.''
Here I use a more explicit construction, getting the
Taylor polynomials of $F$ inductively.
\end{say}

 Let $k$ be a field and $k[[ x_1,\dots,x_n]]$
the ring of formal 
power series over $k$.
For  $f(x_1,\dots,x_n)\in k[[ x_1,\dots,x_n]]$, let
 $$
t_m(f)\in k[[ x_1,\dots,x_n]]/(x_1,\dots,x_n)^{m+1}=:P_n(m)
$$
denote the truncation, mapping a power series to its
degree $m$ Taylor polynomial. We can view
$P_n(m)$ as an affine space over $k$ with natural truncation maps
$t_{m',m}: P_n(m')\to P_n(m)$ for every $m'\geq m$.

The following technical lemma makes it possible to
construct the correct limits of power series.

\begin{lem}\label{good.filter} 
 Let $k$ be a field and $f_i(x_1,\dots,x_n)\in k[[ x_1,\dots,x_n]]$
power series indexed by an infinite set $I$.
There are (nonunique)
infinite subsets $I\supset I_0\supset I_1\supset \cdots$  such that
\begin{enumerate}
\item for every $m$, 
the Zariski closure $Z_m\subset P_n(m)$ of $\{t_m(f_i): i\in I_m\}$ 
is irreducible (over $k$), 
\item for every  Zariski closed $Y\subsetneq Z_m$ there are 
only finitely many $i\in I_m$ such that $t_m(f_i)\in Y$, and
\item for every $m'\geq m$ the truncation maps
$t_{m',m}:Z_{m'}\to Z_m$ are dominant.
\end{enumerate}
\end{lem}

Proof. Apply (\ref{subset.lem}) to 
 $\{t_0(f_i): i\in I\}$ as points in $P_n(0)$
to obtain  $I_0:=I'$.

Assume now that we already have $Z_j\subset P_n(j)$
and $I_j\subset I$ for $j\leq m$ satisfying the properties
(\ref{good.filter}.1--3).

Apply (\ref{subset.lem}) to 
 $\{t_{m+1}(f_i): i\in I_m\}$ as points in $P_n(m+1)$
to obtain  $I_{m+1}:=\bigl(I_m\bigr)'$.
The properties
(\ref{good.filter}.1--2) hold by construction.
The truncation map $t_{m+1,m}:Z_{m+1}\to Z_m$ is defined.
The closure of its image contains all the  points
 $t_{m}(f_i)$ for $ i\in I_{m+1}$. Since
(\ref{good.filter}.2) holds for $Z_m$, we conclude that
 $t_{m+1,m}:Z_{m+1}\to Z_m$ is dominant.\qed

\begin{say}[Generic limits]
 Let $k$ be a field and $f_i(x_1,\dots,x_n)\in k[[ x_1,\dots,x_n]]$
power series indexed by an infinite set $I$.
Let $K\supset k$ be an algebraically closed field
of infinite transcendence degree.
Set $g_{(-1)}=0$ and for $m\geq 1$ 
choose $K$-points
$g_m\in Z_m(K)$ such that
\begin{enumerate}
\item $g_m\in Z_m(K)$ is a generic point 
of $Z_m$ over  $k(g_{m-1})$ (in the sense of Weil, cf.\
 (\ref{weil.gen.pts})), and
\item $g_{m}$ is a lifting of $g_{m-1}$, that is, $t_{m,m-1}(g_{m})=g_{m-1}$.
\end{enumerate}

We can view the $g_m$ as  successive truncations of a
power series $F(x_1,\dots,x_n)\in K[[ x_1,\dots,x_n]]$.
We call any such $F$ a 
{\it generic limit} of the power series $f_i\in k[[ x_1,\dots,x_n]]$.
\end{say}

\begin{thm}\label{limit.series.thm}
 With the above notation, 
$c_0(F)$ is a (Euclidean) limit point of 
$\bigl(c_0(f_i):i\in I\bigr)\subset \r$.
\end{thm}

Proof. By construction, $t_m(F)$ is a generic point of
$Z_m$, hence, as in (\ref{thr.family.say}), there is a Zariski open
$U_m\subset Z_m$ such that
$c_0(h)=c_0\bigl(t_m(F)\bigr)$ for every $h\in U_m$. 
Thus there is a $j\in I_m$ such that
$c_0\bigl(t_m(F)\bigr)=c_0\bigl(t_m(f_j)\bigr)$. Therefore,
using (\ref{taylor.app.cor}) twice, we obtain that 
$$
\Bigl|c_0\bigl(F\bigr)-c_0\bigl(f_j\bigr)\Bigr|\leq 
\Bigl|c_0\bigl(F\bigr)-c_0\bigl(t_m(F)\bigr)\Bigr|+
\Bigl|c_0\bigl(t_m(f_j)\bigr)-c_0\bigl(f_j\bigr)\Bigr|
\leq 2\frac{n}{m+1}.\qed
$$

\begin{lem}\label{subset.lem}
 Let $X$ be a Noetherian topological space and
$\{p_i: i\in I\}$ an infinite  collection of (not necessarily distinct)
  points.
Then there is a (nonunique) subset $I'\subset I$
such that 
\begin{enumerate}
\item the Zariski closure $Z(I')$ of $\{p_i: i\in I'\}$ is irreducible, and
\item for every  Zariski closed $Y\subsetneq Z(I')$
there are only finitely many $i\in I'$ such that $p_i\in Y$.
\end{enumerate}
\end{lem}

Proof. The Zariski closure  of $\{p_i: i\in I\}$
has only finitely many irreducible components.
Pick any, say $X_1$, and set $I_1:=\{i\in I: p_i\in X_1\}$.

If some irreducible closed $X_2\subsetneq X_1$ 
contains infinitely many  $p_i$ for $ i\in I_1$, 
let these be $\{p_i: i\in I_2\}$. We construct 
$X_3 \subsetneq X_2$ similarly, and so on.

By the Noetherian property, eventually we obtain
an infinite subset $I':=I_r\subset I$ such that $X_r$, 
the Zariski closure of $\{p_i: i\in I_r\}$, is irreducible, and
for every  Zariski closed $Y\subsetneq X_r$
there are only finitely many $i\in I_r$ such that $p_i\in Y$.

Note that if a point $p$ appears among the $p_i$
infinitely many times, then $Z(I')=\{p\}$
satisfies the requirements.
\qed

\begin{say}[Generic points \`a la Weil]\label{weil.gen.pts}
 Let $X\subset \a^n_k$ be an irreducible
$k$ variety and $K\supset k$ an algebraically closed extension of
infinite transcendence degree. According to \cite[Sec.IV.1]{weil},
a  {\it generic point} of $X$ is a $K$-point $g_X\in X(K)$
such that a polynomial $p\in k[x_1,\dots,x_n]$ vanishes on
$g_X$ iff it vanishes on $X$. Equivalently,
the restriction map $k(X)\map k(g_X)$ is an isomorphism
where $k(g_X)\subset K$ is the field generated by the
coordinates of $g_X$. 

It is easy to construct generic points as follows.
We may assume that $\dim X=d$ and the projection to the first
$d$ coordinates $\pi:X\to \a^d$ is dominant. Pick
$(p_1,\dots, p_d)\in \a^d$ such that the $p_i\in K$ are
algebraically independent over $k$. Any
$$
(p_1,\dots, p_d, p_{d+1},\dots, p_n)\in \pi^{-1}(p_1,\dots, p_d)\subset X
$$
is a generic point of $X$ over $k$.

Let $f:X\to Y$ be a dominant map of irreducible
$k$-varieties. Let $g_Y\in Y$ be a generic point.
We can view $f^{-1}(g_Y)$ as a $k(g_Y)$-variety.
Any generic point $g_X\in f^{-1}(g_Y)$ as $k(g_Y)$-variety
is also a generic point of $X$ as a $k$-variety.

Thus, given a tower of irreducible $k$-varieties
$$
Z_1\leftarrow Z_2\leftarrow  Z_3\leftarrow \cdots
$$
we can get a compatible system of generic points
$$
g_1\leftarrow g_2 \leftarrow  g_3 \leftarrow \cdots
$$
\end{say}

\section{Taylor polynomials and thresholds}
\label{taylor.sec}

Let $f$ be a holomorphic function and 
assume that $(f=0)$ defines an isolated singularity at the origin.
A result of \cite{hiro} says that there is an $m>0$ such that
if $h$ is any other holomorphic function that agrees with
$f$ up to high order then there is a local biholomorphism 
$\phi:(0\in \c^n)\to (0\in \c^n)$
that takes $f$ to $h$. In particular, the singularities
 $(f=0)$ and $(h=0)$ are analytically isomorphic and 
 all their local analytic invariants
 are the same.

This result completely fails if $(f=0)$ does not have an isolated
singularity. If $p_m$ is a general degree $m$ homogeneous 
polynomial then $(f+p_m=0)$ has only isolated
singularities. In particular, $(f=0)$ and $(f+p_m=0)$
are not isomorphic.

What about their log canonical thresholds? In general, the answer is again
negative. For instance, if $f(z_1,\dots, z_{n-1})$
is any  holomorphic function  then by (\ref{dsum.prop}.2)
$$
c_0\bigl(f(z_1,\dots, z_{n-1})+z_n^m\bigr)=
\min\bigl\{c_0(f)+\tfrac1{m}, 1\bigr\}.
$$
Note  that $(f=0)\subset \c^n$ is very non-isolated; 
it is equisingular along the $z_n$-axis.

The relevant analog of the Hironaka theorem should be  about
holomorphic functions $f$ for which the origin is
isolated ``as far as the log canonical threshold is concerned.'' That is,
 functions $f$ such that
$$
c_0(f)<c_p(f) \qtq{for every $p\neq 0$ near $0$.}
$$

More generally, we
consider the case when $c_0(f)$ is computed by a divisor $E$
whose center on $\hat{\a}^n$ is the origin. 
That is, if there is  a birational morphism $\phi:X\to \hat{\a}^n$ 
and a divisor $E$ on $X$ such that the center of $E$ on
$\hat{\a}^n$  is the origin
and 
$$
c_0(f)=\frac{1+\mult_E \jac(\phi)}{\mult_E f}.
$$
In this case an easy argument  of \cite{df-m}  shows that 
$c_0(f+p)\leq c_0(f)$ if $\mult_0p\gg 1$; see (\ref{upper.anytail.say}).
Here I prove that,
in fact, $c_0(f+p)= c_0(f)$.
For the proof of (\ref{accum.conj})  one needs a more general version,
when certain perturbations of low degree terms are also allowed.
This is considered in Section \ref{accu.sect}.

\begin{thm}\label{anytail.prop}
 Let $f\in k[[x_1,\dots, x_n]]$ be a power series
such that $c_0(f)$ is computed by a divisor $E$
whose center on $\hat{\a}^n$ is the origin.
Let $p\in k[[x_1,\dots, x_n]]$ be a power series
such that
$\mult_0 p> \mult_Ef$.
Then $c_0(f+p)=c_0(f)$.
\end{thm}

The proof  of (\ref{anytail.prop}) 
relies on rather heavy machinery; we need
the full force of the MMP.
 Before starting it, 
let us review a simple argument which gives 
  the inequality  $c_0(f+p)\leq c_0(f)$.

\begin{say}\label{upper.anytail.say}
Take a log resolution  $\pi:X\to \a^n$ and let
$E\subset X$ be a divisor such that
$\pi(E)=0$.
If $\mult_0p>\mult_E  f$
then $ p$ vanishes along $E$ with
multiplicity $>\mult_E  f$, thus
$$
\mult_E  (f+p)=\mult_E  f.
\eqno{(\ref{upper.anytail.say}.1)}
$$
If we choose $E$ such that
$$
c_0(f)=\frac{1+\mult_E \jac(\pi)}{\mult_E  f},
$$
then we obtain that 
$$
c_0(f+p)\leq \frac{1+\mult_E \jac(\pi)}{\mult_E  (f+p)}=
\frac{1+\mult_E \jac(\pi)}{\mult_E  f}=c_0(f).
\eqno{(\ref{upper.anytail.say}.2)}
$$

In this argument, the dependence on $f$ is rather
subtle. We need a complete log resolution
in order to choose the right divisor $E$. 
Unfortunately, $X\to \a^n$ is not a log resolution
for the perturbed hypersurface $(f+p=0)$, except 
in the isolated singularity case.
A priori, the log canonical threshold of $f+p$ may  be computed
by a divisor $E'$ which does not even appear
on $X$. The values of $\mult_{E'}\jac(\pi)$ and $\mult_{E'}  f$
may be completely different from 
 $\mult_{E}\jac(\pi)$ and $\mult_{E}  f$.
Thus we can not just write $f=(f+p)+(-p)$ and
obtain the reverse inequality
$c_0(f)=c_0\bigl((f+p)+(-p)\bigr)\leq c_0(f+p)$.

Note, however, that we are on the right track.
By the ACC conjecture (\ref{acc.conj}), there are no
log canonical thresholds
 in some interval  $\bigl(c_0(f)-\epsilon, c_0(f)\bigr)$, and,
by (\ref{taylor.app.cor}), $c_0(f+p)\geq c_0(f)-{n}/{\mult_0p}$.
Thus  $ c_0(f+p)\geq c_0(f)$ if  $\mult_0p>n/\epsilon$.
\end{say}

\begin{rem} Assume that  $c_0(f)$ is computed by a divisor $E$
whose center $Z(E)$ on $\hat{\a}^n$ is not the origin.
If $p$ vanishes along $Z(E)$ with multiplicity $>\mult_Ef$,
then $c_0(f+p)=c_0(f)$ should  hold.
The problem with the proof is that usually $E$ can not
be realized on an algebraic variety $X\to \a^n$.

 It should be possible to obtain 
(\ref{mmp.lem}) using the MMP
and a log resolution of $(\hat{\a}^n,(f=0))$. However,
 the standard references seem to
assume that we consider MMP over a base which itself
is a variety or an analytic space.
\end{rem}

\section{Proofs}\label{accu.sect}

 For the proofs it is more convenient to
change to additive notation. If $X$ is a complex manifold
and $D_i\subset X$ are divisors with local equations $f_i$
then we say that $(X, \sum c_iD_i)$ is {\it lc} (or {\it log canonical})
if $\bigl(\prod_i |f_i|^{-c_i}\bigr)^{1-\epsilon}$ is
locally $L^2$ for every $\epsilon>0$. 

In the sequel we also need these concepts when $X$ itself is
singular. See \cite[sec.2.3]{km-book} for a good introduction.

\begin{say}[Proof of (\ref{anytail.prop})]
\label{pf.of.anytail.prop}
 By  (\ref{mmp.lem}), there is a proper birational morphism
$\pi:X\to {\a}^n$ such that $X$ is $\q$-factorial and
$E$ is (birational to) the unique exceptional divisor of $\pi$.
Completing over the origin, we obtain
a proper birational morphism
$\hat{\pi}:\hat X\to \hat{\a}^n$ such that 
$E$ is (birational to) the unique exceptional divisor $\hat E$ of $\hat\pi$
and $\hat E$ is $\q$-Cartier. Thus
$\bigl(\hat X, \hat E+c\cdot \hat\pi^{-1}_*(f=0)\bigr)$ is lc
where $c=c_0(f)$ and $\hat\pi^{-1}_*$ denotes the birational
transform of a divisor.

We need to prove that 
$\bigl(\hat X, \hat E+c\cdot \hat{\pi}^{-1}_*(f+p=0)\bigr)$ is also lc. 
As a first step, we claim that
$$
\bigl(\hat{\pi}^{-1}_*(f=0)\bigr)|_{\hat E}=
\bigl(\hat{\pi}^{-1}_*(f+p=0)\bigr)|_{\hat E}.
\eqno{(\ref{anytail.prop}.1)}
$$
In order to prove this, let us remember
how  we compute $\hat{\pi}^{-1}_*(f=0)$.
The zero divisor of $f\circ \hat{\pi}$ is
$r\hat E+\hat{\pi}^{-1}_*(f=0)$ for some $r>0$. 
 It is better to work with Cartier divisors,
so choose $m>0$ such that
$mr\hat E$ and $m\cdot \hat{\pi}^{-1}_*(f=0)$ are both Cartier.
Thus, if $y=0$ is a local equation of $mr\hat E$ then
the local equation of  $m\cdot \hat{\pi}^{-1}_*(f=0)$ is
$\bigl(f^m\circ \hat{\pi}\bigr)/y=0$. 
By assumption, $\mult_{\hat E}p\geq \mult_0p>
\mult_{\hat E}f$,
thus every term of 
$$
(f+p)^m\circ \hat{\pi}=
\sum_i \binom{m}{i}\bigl(f^i\cdot p^{n-i}\bigr)\circ \hat{\pi}
$$
  is divisible by $y$ and
all except $\bigl(f^m\circ \hat{\pi}\bigr)/y$ vanish along $\hat E$.
Thus (\ref{anytail.prop}.1) holds.

By the precise version of inversion of adjunction  
(\ref{prec.inv.adj}) this means that 
$$
\bigl(\hat X, \hat E+c\cdot \hat\pi^{-1}_*(f=0)\bigr)\qtq{is lc iff}
\bigl(\hat X, \hat E+c\cdot \hat{\pi}^{-1}_*(f+p=0)\bigr)\qtq{is.}
$$ 
Again we have to be careful with the formal case.
The easiest is to notice that
$$
\bigl(\hat{\pi}^{-1}_*(f=0)\bigr)|_{\hat E}=
\bigl({\pi}^{-1}_*(f=0)\bigr)|_{E},
$$
thus the algebraic version of inversion of adjunction  
gives that 
$$
\bigl( X,  E+c\cdot \pi^{-1}_*(t_m(f+p)=0)\bigr)\qtq{is lc for $m>\mult_Ef$.}
$$
This implies that $\bigl(\a^n, c\cdot (t_m(f+p)=0)\bigr)$
is lc, thus $c_0\bigl(t_m(f+p)\bigr)\geq c$ for $m\gg 1$.

By (\ref{taylor.app.cor}) this implies that
$$
c_0(f+p)=\lim_{m\to\infty} c_0\bigl(t_m(f+p)\bigr) \geq c=c_0(f). \qed 
$$
\end{say}

\begin{lem}\label{mmp.lem}
 Let $\hat D\subset \hat{\a}^n$ be a divisor and
$E$ a divisor over $\hat{\a}^n$ such that the center of $E$
on $\hat{\a}^n$ is the origin and $E$ computes the
log canonical threshold of $\hat D$.

Then there is a proper birational morphism $\pi:X\to \a^n$
with only 1 exceptional divisor, which is (birational to) $E$.
\end{lem}

Proof. The first step is to note that $E$ is an algebraic divisor.
That is, there is a proper birational morphism $g:Y\to \a^n$
such that $g$ is an isomorphism outside the origin, $E$ is a divisor
on $Y$ and $\ex(g)$ is a snc divisor.
This follows from the resolution of indeterminacies of  maps.
See \cite[p.113]{ksc} for an elementary proof.
 
If $\hat D$ is defined by the power series $f$,
set $D_m:=(t_mf=0)$. Set $c=c_0(f)$.
Write
$$
c\cdot g^* \hat D=K_{\hat Y}+\Delta+c\cdot g^{-1}_*\hat D
\qtq{and}   
c\cdot g^* D_m=K_Y+\Delta_m+c\cdot g^{-1}_*D_m.
$$
Note that   $\Delta_m=\Delta$ for $m\gg 1$,
thus
$E$ appears in $\Delta_m$ with coefficient 1.
The pair $(\a^n, c\cdot D_m)$ is not yet known to be
lc, but, by (\ref{taylor.app.cor}),  $c_0(D_m)$ converges to $c=c_0(\hat D)$.
Thus, for  $c'<c$  and $m\gg 1$,
$(\a^n, c'\cdot D_m)$ is klt  and if we write
 $$
c'\cdot g^* D_m=K_Y+\Delta'_m+c'\cdot g^{-1}_*D_m
$$
then $E$ appears in $\Delta'_m$ with coefficient $>0$.
Apply 
(\ref{1-div.lem}) to $E$ and $(\a^n, c'\cdot D_m)$
to get $\pi:X\to \a^n$ as required.\qed

\begin{lem}\label{1-div.lem} Let $X$ be a variety and $D$ a
$\q$-divisor on $X$  such that $(X,D)$ is klt.
Let $E$ be a divisor over $X$ such that
$0\leq a(E,X,D)<1$. 
Then there is a proper birational morphism $\pi:X_E\to X$
with only 1 exceptional divisor, which is (birational to) $E$.
\end{lem}

Proof. Let $g:Y\to X$ be a log resolution such that
$E$ is a divisor on $Y$. 
 Write
 $$
K_Y+eE+A-B+ g^{-1}_*D\simq g^*\bigl(K_X+D),
$$
where $e=a(E,X,D)$, the $\q$-divisors  $A,B$ are effective and have
no common components and do not contain $E$.
 For some $0<\eta\ll 1$, run the
$(Y, eE+(1+\eta)A)$-MMP \cite{bchm}
to obtain $\pi:X_E\to X$.
Note that
$$
K_Y+eE+(1+\eta)A\simq \eta A+B+g^*\bigl(K_X+D)
$$
and on the right hand side $\eta A+B$ contains every $g$-exceptional divisor
with positive coefficient, save $E$. 
Thus the restriction of any birational transform
of $K_Y+eE+(1+\eta)A$ to the (birational transform) of
$E$ is $\q$-linearly equivalent to an
effective divisor plus a pulled-back divisor, hence we never contact $E$.
On the other hand, an effective exceptional divisor
is never relatively nef, thus we have to contract
$\supp (\eta A+B)$.
Thus $\pi:X_E\to X$
has only one exceptional divisor, which is (birational to) $E$.
\qed

\begin{say}[Precise  inversion of adjunction]\label{prec.inv.adj}

As stated in (\ref{ot-inv}), inversion of adjunction is only an
inequality. It is possible to make it into an equality.
Assume that $X$ is smooth,  $E\subset X$ is a hypersurface
and $\Delta$ an effective $\q$-divisor which does not contain $E$.
 The precise  inversion of adjunction says that
$$
(X, E+\Delta) \qtq{is lc near $E$} \Leftrightarrow
(E, \Delta|_E) \qtq{is lc.}
\eqno{(\ref{prec.inv.adj}.1)}
$$
In (\ref{pf.of.anytail.prop}) we use a version where $X$ is singular.
For the  statements and proofs see
\cite[Sec.17]{f&a} and \cite{kawakita}.
\end{say}

\begin{say}[Proof of (\ref{main.accu.thm})] 
It follows from (\ref{dsum.prop}.2) that the set of accumulation points of
${\mathcal{HT}}_n$ contains
 ${\mathcal{HT}}_{n-1}\setminus\{1\}$.
It is thus enough to prove that it is contained in
${\mathcal{HT}}_{n-1}$.

Let $f_i\in k[[x_1,\dots,x_n]]$ be power series
such that $c_0(f_i)$ is a nonconstant sequence
converging to some $c\in \r$. 
By passing to a subsequence we may assume that
$c_0(f_i)\neq c_0(f_j)$ for $i\neq j$.
By
(\ref{limit.series.thm}), we get a  power series $F\in K[[x_1,\dots,x_n]]$
such that $c_0(F)=c$.

If $c_0(F)$ is computed by a divisor whose center $Z(E)$ is not the origin,
then localizing at the generic point of $Z(E)$ and completing
gives a complete, regular, local ring of dimension
$n-\dim Z(E)$ and a power series $F^*$ such that
$c_0(F^*)=c_0(F)$. Thus
$c_0(F)\in {\mathcal{HT}}_{n-1}$ and we are done.
Otherwise 
$c_0(F)$ is computed by a divisor $E$ whose center is the origin.
By (\ref{semicont.thm}) this implies that
  $c_0(F)=c_0(f_j)$ for some infinite subsequece $i_1<i_2<\cdots$,
a contradiction. \qed
\end{say}

\begin{prop}\label{semicont.thm}
  Let $K\supset k$ be a field extension
and $F(x_1,\dots,x_n)\in K[[ x_1,\dots,x_n]]$ a power series.
Assume that $c_0(F)$ is computed by a divisor $E$
whose center on $\hat{\a}^n$ is the origin. 
 Let $Z_m\subset P_n(m)$ denote the
 $k$-Zariski closure of $t_m(F)$.
Then there is an $m\geq 0$ and a nonempty open subset
$U_m\subset Z_m$
such that if $f(x_1,\dots,x_n)\in K[[ x_1,\dots,x_n]]$ is any power series
such that $t_m(f)\in U_m$
then $c_0(f)=c_0(F)$.
\end{prop}

Proof. By (\ref{mmp.lem}), there is a proper birational morphism
 $\phi:X\to {\a}^n$
defined over $K$  with a unique $\phi$-exceptional divisor,
which is (birational to) 
$E$.

The data $\phi, X$ and  $E$
 are defined over a finitely generated subextension
of $K/k$, hence over $k\bigl(t_m(F)\bigr)$ for all $m\gg m_K$ for some $m_K$.
Set $m_E:=\mult_E (F\circ \phi)$ and choose any $m> \max\{m_E,m_K\}$. 

Since $\phi, X$ and  $E$
 are defined over the generic point of $Z_m$,
 there is a 
Zariski open subset  $U_m\subset Z_m$
such that $\phi, X$ and  $E$
can be extended to be defined over $U_m$.
Moreover, we may assume that for any $u\in U_m$, the
resulting
$\phi(u):X(u)\to {\a}^n$ is birational,
$E(u)\subset X(u)$ is a divisor with the same
discrepancy as $E\subset X$,
  $F(u)(\phi_1(u),\dots, \phi_n(u))$
vanishes along $E(u)$ with multiplicity $m_E$.

By shrinking $U_m$ if necessary, we may also assume
that if $t_m(f)\in U_m$
then  $c_0\bigl(t_m(f)\bigr)=c_0\bigl(t_m(F)\bigr)$.

Let now $f$ be any power series
with $t_m(f)\in U_m$ and apply
(\ref{anytail.prop}) 
first to $t_m(F)$ and $p:=F-t_m(F)$ and then 
to  $t_m(f)$ and $p:=f-t_m(f)$. 
We obtain that
$$
\begin{array}{rcccl}
c_0(F)&=&c_0\bigl(t_m(F)+(F-t_m(F))\bigr)&=&c_0\bigl(t_m(F)\bigr),\qtq{and}\\
c_0(f)&=&c_0\bigl(t_m(f)+(f-t_m(f))\bigr)&=&c_0\bigl(t_m(f)\bigr).
\end{array}
$$
Since $c_0\bigl(t_m(f)\bigr)=c_0\bigl(t_m(F)\bigr)$
by the choice of $U_m$, we conclude that
$c_0(f)=c_0(F)$.\qed

\section{Technical comments on inversion of adjunction}
\label{techn.sec}

\begin{say}[Inversion of adjunction: the isolated singularity case]

Let us consider inversion of adjunction (\ref{ot-inv}) in case
$D:=(f=0)$ has an isolated singularity. 
Then there is an algebraic hypersurface $D'$ with an isolated
sigularity at the origin such that 
a formal change of coordinates transforms $D$ to $D'$.
Let $\pi:X\to \a^n$ be an algebraic resolution
of $(\a^n, D')$ (in a neighborhood of  the origin) 
that is an isomorphism outside the origin.
By completion, we get a log resolution
$\hat{\pi}:\hat{X}\to \hat{\a}^n$   
of $(\hat{\a}^n, D)$ which is an isomorphism outside the origin.
Up to coordinate change we may assume that
$\hat{L}\subset \hat{\a}^n$ is a coordinate subspace
that is the completion of a linear subspace $L\subset \a^n$.

Set $D_m:=(t_m(f)=0)$. Then $\pi:X\to {\a}^n$
is a log resolution
of $({\a}^n, D_m)$ in a neighborhood of the origin for
$m\gg 1$.

The proof of inversion of adjunction given in \cite[Sec.17]{f&a}
shows that for each $m\gg 1$ there is an irreducible component
 $E_m\subset \pi^{-1}_*L\cap \pi^{-1}(0)$
such that 
$a(E_m, L, c_0(f)\cdot D_m|_L)\geq 1$.
Since $\pi^{-1}_*L\cap \pi^{-1}(0)$ has only finitely many
irreducible components, by passing to a subsequence we may assume that
$E=E_m$ does not depend on $m$.
By (\ref{upper.anytail.say}.1), for divisors with center at the origin
the discrepancy stabilizes, hence
$$
a(E, L, c_0(f)\cdot D|_L)=a(E, L, c_0(f)\cdot D_m|_L)
\qtq{for $m\gg 1$.}
$$
Thus $c_0(f|_L)\leq c_0(f)$.\qed
\end{say}

As in (\ref{dk.proof}), this implies the following
special case of (\ref{dk.thm}):

\begin{lem} Let $f,g\in k[[x_1,\dots,x_n]]$ be formal power series
with  isolated singularities at the origin. Then
$c_0(f+g)\leq c_0(f)+c_0(g)$. \qed
\end{lem}

\begin{say}
Now we are ready to prove the theorems of Section \ref{basic.props.sec}
 in the formal case.
As noted in (\ref{formal.say}), we only need to show that
$$
\liminf_{m\to\infty} c_0\bigl(t_m(f)\bigr)\geq c_0(f).
$$
Let $f\in k[[x_1,\dots,x_n]]$ be a formal power series
with degree $m$ Taylor polynomial $t_m(f)$.
Let $h$ be  a general degree $m+1$ homogeneous polynomial.
Then both
$$
t_m(f)+h\qtq{and}  \bigl(f-t_m(f)\bigr)-h
$$
have isolated singularities at the origin and so
$$
c_0(f)\leq c_0\bigl(t_m(f)+h\bigr)+c_0\bigl(f-t_m(f)-h\bigr).
$$
On the other hand, by the algebraic version of
(\ref{dk.thm}), 
$c_0\bigl(t_m(f)+h\bigr)\leq c_0\bigl(t_m(f)\bigr)+c_0(h)$. Therefore,
using (\ref{compute.lcth}.5) we get that 
$$
c_0(f)\leq c_0\bigl(t_m(f)\bigr)+c_0(h)+c_0\bigl(f-t_m(f)-h\bigr)
\leq c_0\bigl(t_m(f)\bigr)+\tfrac{n}{m+1}+\tfrac{n}{m+1}. \qed
$$
\end{say}

 \begin{ack}  I thank R.\ Lazarsfeld and M.\ Musta{\cedilla{t}}{\v{a}}
for useful e-mails and corrections. 
Partial financial support  was provided by  the NSF under grant number 
DMS-0500198. 
\end{ack}

\bibliography{refs}

\vskip1cm

\noindent Princeton University, Princeton NJ 08544-1000

\begin{verbatim}kollar@math.princeton.edu\end{verbatim}

\end{document}